\documentclass{article}

\usepackage{arxiv}

\usepackage[utf8]{inputenc} 
\usepackage[T1]{fontenc}    
\usepackage{hyperref}       
\usepackage{url}            
\usepackage{booktabs}       
\usepackage{amsfonts}       
\usepackage{nicefrac}       
\usepackage{microtype}      
\usepackage{lipsum}
\usepackage{graphicx}
\graphicspath{ {./images/} }

\title{Optimal Competition Resolution Rule for Buslaev Controlled Binary Chain
}

\author{
  Tatashev A.G. \\
  Department of Higher Mathematics\\
  Moscow Automobile and Road Construction\\
  State Technical University (MADI) \\
  Moscow, Leningradsky avenue, 64, Russia  \\
and MTUCI, Moscow, Aviamotornaya str., 8-a, Rissia\\  
  \texttt{a-tatashev@yandex.ru} \\
   \And
 Yashina M.V. \\
  Department of Higher Mathematics\\
  Moscow Automobile and Road Construction\\
  State Technical University (MADI) \\
  Moscow, Leningradsky avenue, 64, Russia  \\
and MTUCI, Moscow, Aviamotornaya str., 8-a, Rissia\\
  \texttt{mv.yashina@madi.ru} \\
}

\begin{document}
\maketitle
\begin{abstract}
A dynamical system, called a  binary closed chain of contours, is studied. The dynamica system belongs to the class of  Buslaev networks. 
The system contains  $N$  {\it contours.} There two cells and a  particle in each contour. 
There two adjacent contours for each contour. There is a common point of  adjacent contours. This common point is called a node. The node is located between the cells. In the deterministic version of the system, at any discrete moment, each  particles moves to the other cell of the contour if there is no delay. The delays are due to that two particles may not pass through the common node 
simultaneously.  If two particles try to cross the same node, then a {\it competition} 
occurs, and only one of these particles moves  in accordance with a prescribed competition resolution rule.  In the stochastic version of the system, each particle moves with the probability $1-\varepsilon,$ if the system is in the state such that, in the same state of the deterministic system, this particle moves.   where $\varepsilon$ is a small value. We have obtained a competition resolution rule such 
that the system results in a state such that all particles move without delays in present time and in the future (a state of free movement), and the system results in the state of free movement over a minimum time. The expectation of the number of the $i$th particle transitions per a time unit is called the {\it average velocity of this particle,}  $v_i,$ $i=1,\dots,N.$ For the stochastic version of 
the system, under the assumption that $N=3,$ we have proved the following. For the optimal rule, the average velocity of particles is equal to $v_1=v_2=1-2\varepsilon+o(\varepsilon)$  $(\varepsilon\to 0).$ For the left-priority rule, which is studied earlier, the average velocity of particles equals $v=v_1=v_2=\frac{6}{7}+o(\sqrt{\varepsilon}).$

\end{abstract}

\keywords{
Dynamical systems \and  Buslaev networks \and Binary chain \and 
Optimal control \and Competition resolution rule 
}

\section{Introduction}

\section*{Introduction.} A class of mathematical models is formed by dynamical systems such that in these systems particles move on closed or infinite one-dimensional lattices. A class of mathematical traffic models is formed by dynamical systems such that, in these system, particles move on closed or infinite one-dimensional lattices. For simple versions of these systems, results have been obtained, for example in [1]--[8]. More complicated systems of these class, in particularly, models with network structure, are investigated mainly by simulation. Models of these class can be interpreted in terms of cellular automata [9], asynchronous exclusion processes [10 ],  or synchronous exclusion processes, [9]. In [11]--[12], traffic models are investigated such that, in these models, particles move in two perpendicular directions on a toroidal lattice. 
 
In [13], a dynamical system is studied such that particles move in channels of this system under a prescribed plan.       

In [14], the concept of cluster movement for mathematical models of traffic were introduced. In the discrete versions, the clusters are groups of particles located in adjacent cells and moving simultaneously. In the continuous versions the clusters are moving segments.  A.P.~Buslaev has introduced a class of dynamical systems called contour networks (Buslaev networks) [15]. A contour networks 
contains contours. There are common points of adjacent contours called nodes. There are particles or clusters on the contours. The particles (clusters) move in accordance with a prescribed rule. Delays occur at the nodes.  The delays are due to that the particles (clusters) may not pass through the node simultaneously. Contour networks may be traffic models and may have other applications. In particular, contour networks also used for modeling the work of communication systems.  In [15]--[23] analytical results were obtained for contour networks. 

In [15], a dynamical system, called a closed binary chain of contours, is considered. There are two cells and a particle in each contour. There are two nodes in each contour. These nodes are common points for this node and two adjacent contours. If two particles try to pass through the same node simultaneously, then a competition occurs. The paper [16] studies binary chains with three rules for competition resolution. These rules are a stochastic rule in accordance with each of two competing particles passes through the node first (wins the competition) equiprobably, the left-priority rule (the particle, located on the left wins), and the rule such that, in accordance with this rule the particle, located on the contour with an even index, wins.   
The paper [17] studies a stochastic version of the system consider in [16]. If the stochastic system is in the state such that a particle moves in deterministic version, then the probability that the particle moves is $1-\varepsilon.$  The behavior of the system is studied under the assumption that $\varepsilon$ tends to 0. A formula for the average velocity of particles has been obtained for the competition resolution rule in accordance the particles located in the contours with even indices win the competitions. In [18], a closed binary chain of contours is considered such that,
in this system, no competing particle moves in the current moment and in the future. In [18], a generalization of the binary chain of contours is also considered. Namely, there are $m$ cells and a particle is in each contour, where $m$ is an even number. In [19], a closed chain of contour is considered. There are $m$ cells and a cluster, containing $l$ particles, is considered, where $m$ is an even number, $l<m.$  A continuous counterpart of this system is considered in [20].    

This paper consider a binary chain of contours. The optimal competition resolution rule has been found. Namely, the system with this rule results in a state of free movement from any initial state over minimum time interval. This rule is called the long cluster rule. A stochastic version of this system is also considered. In this version, the probability that an attempt of a particle to move is realized is $1-\varepsilon.$ The behavior of the system is studied under the assumption that $\varepsilon$ tends to 0.   
 
 \section*{2. Description of the system}

\hskip 18pt We consider system containing $N$ contours, Fig. 1. The indices of the contours are $0,1,\dots,N-1.$ There are two cells on each contour. These cells are the lower cell~0 and the upper cell~1. We say that the contour $i$ is in the state~$j$ if the particle of this contour (particle~$i)$ is in the cell~$j,$  $j=1,2.$ For each contour, there are two adjacent contours. One of these adjacent contours is located on the left, and the other adjacent contour is on the right. The contour $i-1$ (subtraction by modulo $N)$ is located to the left of the contour $i,$ and the contour $i+1$ (addition by modulo $N)$ is located to the right of the contour $i,$ $i=0,1,\dots,N-1.$ There is a common point of the nodes $i,$ $i+1$ called the node $(i,i+1)$ (addition by modulo $N).$ Passing from the cell~0 to the cell~1, the particle of the contour $i$ (the particle $i)$ crosses the node $(i,i+1),$ and, passing from the cell~0 to the cell~1, the particle  $i$ crosses the cell $(i-1,i)$ (counter-clockwise movement), $i=0,1,\dots,N-1.$

\begin{figure}[ht!]
\centerline{\includegraphics[width=400pt]{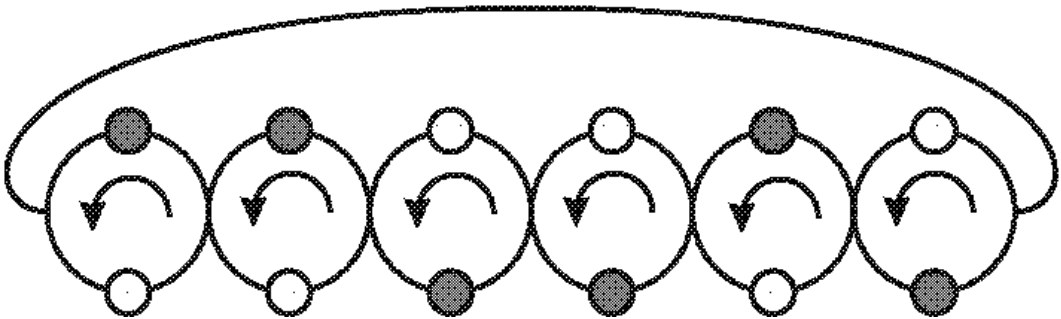}}
\caption{ Binary closed chain of contours}
\end{figure}

In the deterministic version, at any discrete moment, each particle passes to the other cell of the contour if there is no delay. If two particles try to pass through the same node, then a      
{\it competition} occurs. In this case, only one of the competing particles moves chosen in accordance with a prescribed rule for competition resolution.  

A state of the system is a cyclic vector
$$(d_0,d_1,\dots,d_{N-1}),$$  
where $d_i$ is the index of the cell in which the particle $i$ is located, $i=0,1,\dots,N-1.$ A set of `'1''s, located in the adjacent positions of the cyclic state is called a {\it 1-cluster).} The clusters of `'1''s are separated by `'0''s. 
The number of `'1'' in a 1-cluster is called the length of this 1-cluster. A {\it 0-cluster} and its length is defined analogously.   

The initial state of the system is prescribed.

\section*{3. Description of the stochastic system}
  
\hskip 18pt  In the stochastic version, the probability that a particle moves is $1-\varepsilon$ if the system is in a state such that, in the deterministic system, the particle moves. If, in the corresponding state of the deterministic system, a particle does not move, the particle does not move in the stochastic system. 

\section*{4. Description of the long cluster rule} 

\hskip 18pt Let $l_0(x)$ be the length of the longest 0-cluster under the system is in the state $x$ and $l_1(x)$ be the length of the longest 1-cluster under this assumption.    If the system is in the state $x$ and $l_0(x)>l_1(x),$ then, in accordance with the {\it long cluster rule,} for any pair of competing pair, the particle located on the left, i.e., the particle, passing from the cell 0 to the cell 1, moves.  If $l_0(x)<l_1(x),$ then, in accordance with the {\it long cluster rule,} for any pair of competing pair, the particle located on the right, i.e., the particle, passing from the cell 1 to the cell 0, moves. If $l_0(x)=l_1(x),$ then we suppose that the competing particle located on the left moves.

In the system with the long cluster rule either, at any step, the length of the longest 0-cluster is decreased by 1 or the the length of the longest 1-cluster is decreased by 1 

\section*{5. Concept of the average velocity}

\hskip 18pt Let $H(t)$ be the expectation of the total number of transitions of particles in the time interval $(0,t).$ The limit 
$$v_i=\lim\limits_{t \to \infty}\frac{H(t)}{t},\ i=0,1,\dots,N-1,$$
is called the {\it average velocity of the particle} $i$, $i=0,1,\dots,N-1.$

The value 
$$v=\frac{1}{N}\sum_{i=0}^{N-1}  v_i$$
is called the {\it average velocity of particles.}

We say that the system is in a state of free movement at time~$t_0$ if all particles move at any time $t\ge t_0.$ If the system results in a state of free movement, then $v_0=v_1=\dots=v_{N-1}.$

\section*{6. Optimization property of the long cluster rule}

\hskip 18pt  Denote by $S_0$ the long cluster rule.     
\vskip 3pt
Suppose
$$l(x(t))=\min(l_0(x(t)),l_1(x(t))),$$  
where $x(t)$ is the state of the system at time $t.$

\vskip 3pt
{\bf Lemma 1. } {\it For any competition competition resolution rule,
$$l(x(t+1)\ge l(x(t))-1.\eqno(1)$$
\vskip 3pt
Proof.}
 Suppose $l(x(t))\ge 1,$ and, in the vector $x(t),$ the longest 0-cluster occupies the positions $i_0,1_0+1,\dots,i_0+s-1$ 
(addition by modulo~$N$).

Then, at time $t,$ 
the particle $i_0+s-1$ 
is a  competing  particle and the particles $i_0,1_0+1,\dots,i_0+l_0-1$ pass to the cells with the index~1.

 Thus there is a 1-cluster of length not less than $l_0(x(t))-1$ in the vector $x(t+1).$ Therefore, 
$$l_1(x(t+1))\ge l_0(x(t))-1.\eqno(2)$$     
The proof of the inequality 
$$l_0(x(t+1))\ge l_1(x(t))-1\eqno(3)$$
is analogous. Using (2), (3) we get (1). Lemma 1 has been proved.

\vskip 3pt
{\bf Lemma 2.} {\it Suppose $l(x(t))\ge 1.$ Then, for the long cluster rule $S_0$ and any state $x,$ the following equality is true
$$l(x(t+1))= l(x(t)).$$ 
\vskip 3pt
Proof.} Suppose that, in the cyclic vector $x(t)$ the longest 0-cluster is in the positions $i_0,1_0+1,\dots,i_0+s-1$ (addition by modulo $N).$ In accordance with the rule $S_0,$ at time $t+1$ the particles $i_0-1,$ $i_0+s-1$ are in the cells~$0$ and the particles $i_0,1_0+1,\dots,i_0+l_0-2$ are in the cells~1. From this and Lemma~1, Lemma~2 follows. 
\vskip 3pt

Assume  that $a(x,S)$ is the number of steps such that over these steps the system results in a state of free movement for the initial state $x$ and the competition resolution rule $S.$ If the system with the competition resolution rule $S$ does not results in a state of free movement, then we suppose that $a(x,S)=\infty.$ 

\vskip 3pt
{\bf Theorem 1.} {\it For any competition resolution rule $S$ and initial state $x$ the inequality  
$$a(x,S)=l(x)\le a(x,S_0)$$ 
is true, i.e., the long cluster rule minimizes the time interval over that the system results to a state of free movement from any initial state.}
\vskip 3pt
Theorem 1 follows from Lemmas 1 and 2.

\section*{7. Stochastic system with three contours and left-priority rule}

\hskip 18pt Suppose that $N=3.$ There are 8 states
$$E_0=(0,0,0),\  E_1=(0,0,1),\  E_2=(0,1,0),\ E_3(0,1,1),$$
$$E_4=(1,0,0),\  E_5=(1,0,1),\  E_6=(1,1,0),\ E_7(1,1,1),$$

Assume that the state space is divided into three subsets  
$$G_1=\{E_0,\},$$
$$G_2=\{E_7,\}$$
$$G_3=\{E_3,E_5,E_6,\},$$
$$G_4=\{E_1,E_2,E_4.\}$$
Due to symmetry, the probability of the transition from the set $G_i$ to the state $G_j$ does not depend on the system state belonging to the state $G_i,$ $i,j=1,2,3,4.$  Thus $G_0,G_1,G_2,G_3$ are states of a Markov chain (macrostates). Thus the process of the system work is a Markov chain such that the states of this chain form a unique non-periodic communicating class, and therefore there exist positive stationary probability, and these probabilities do not depend on the initial state, [24]. 

Denote by $p_{ij}$ the probability that the system passes from the macrostate $G_i$ to the set  $G_j$ over one step, $i,j=1,2,3.$ We have
$$p_{11}=o(\varepsilon),\ p_{12}=1-3\varepsilon+o(\varepsilon),\   p_{13}=3\varepsilon+o(\varepsilon),\  p_{14}=o(\varepsilon),\eqno(1)$$
$$p_{21}=1-3\varepsilon+o(\varepsilon),\   p_{22}=o(\varepsilon),\  p_{23}=o(\varepsilon),\  p_{24}=3\varepsilon+o(\varepsilon),\eqno(2)$$ 
$$p_{31}=0,\ p_{32}=\varepsilon+o(\varepsilon),\ p_{33}=1-2\varepsilon+o(\varepsilon),\  p_{34}=\varepsilon+o(\varepsilon),\eqno(3)$$
$$p_{41}=0,\   p_{42}=1-2\varepsilon+o(\varepsilon),\  p_{43}=2\varepsilon_o(\varepsilon),\ p_{44}=o(\varepsilon)\eqno(4).$$
The states of the Markov chain form a unique non-periodic communicating class and, hence, [30],
there exist stationary probabilities of all states independent of the initial state. Denote by $p_i$ the stationary probabiility of the macrostate $G_i,$ $i=1,2,3,4.$ Using (1)--(4), we get equations for the stationary probabilities of the macrostates
$$p_1=(1-3\varepsilon)p_2+o(\varepsilon),\eqno(5)$$
$$p_2=(1-3\varepsilon)p_1+\varepsilon p_3+(1-2\varepsilon)p_4+o(\varepsilon),\eqno(6)$$
$$p_3=3\varepsilon p_1+(1-2\varepsilon)p_3+2\varepsilon p_4+o(\varepsilon),\eqno(7)$$
$$p_4=3\varepsilon p_2+\varepsilon p_3+o(\varepsilon),\eqno(8).$$  
$$p_1+p_2+p_3+p+4=1.\eqno(9)$$
Using (5)--(9), we get
$$p_1=\frac{2}{7}+o(\sqrt{\varepsilon}),\eqno(10)$$
$$p_2=\frac{2}{7}+o(\sqrt{\varepsilon}),\eqno(11)$$
$$p_3=\frac{3}{7}+o(\sqrt{\varepsilon}),\eqno(12)$$
$$p_4=o(\sqrt{\varepsilon}).\eqno(13).$$
If the system is in the macrostate $G_1$ or $G_2,$ then all partiicles move. If the system is in the macrostate $G_3,$ then a delay of one paricle occurs. From this, using (10)--(13), we get the following proposition. 
\vskip 3pt
{\bf Proposition 1.} {\it If $N=3$ and the left-priority competition resolution rule is prescribed, then the average velocity of any particle is equal to
$$v=v_1=v_2=v_3=\frac{6}{7}+o(\sqrt{\varepsilon}).$$ 
}
\vskip 3pt
{\bf Remark 1.} {\it Suppose that $N$ is an even number, and the competition resolution rule the particles located in the contours with even numbers win the competitions 
({\it odd-even rule}). Then, in accordance with results, obtained in [21], the average velocities of particles are equal to
$$v_i=1-\varepsilon,\ i=0,2,\dots,N-2, $$      
$$v_i=\frac{3}{4}-o(\sqrt{\varepsilon}),\  i=1,3,\dots,N-1,$$
$$v=\frac{7}{8}-o(\sqrt{\varepsilon}).$$
}

\section*{8. Stochastic system with rule of long cluster}

\hskip18pt We consider a stochastic closed chain of contours with the long cluster rule. 

Suppose that the number of contours equals $N,$ and the probability that the particles moves with the particle moves, if the system is in the state such that the particle may move, is equal to $1-\varepsilon.$  
\vskip 3pt
{\bf Theorem 2.} {\it For a stochastic closed chain, the average velocity of clusters is equal to 
$$v=1-\varepsilon+o(\varepsilon),\ \varepsilon\to 0.$$
\vskip 3pt
Proof.} From the state $(0,\dots,0),$ the system results in any prescribed state with a positive probability, in particular, the system may stay in the state $(0,\dots,0).$ From any state, with a positive probability, the system results in the sttate $(1,\dots,1),$  and, from the state $(1,\dots,1),$ the system results in the state $(0,\dots,0).$ 

Assume that the system state space is divided into sets $S_0,S_1,\dots,S_{[N/2]},$ where $([N/2]$ is the integral part of the number $N/2).$ The set $S_i$ is the set of all $x$ such that $l(x)=i,$ $i=0,1,\dots,[N/2].$ Denote by $P_i$ the stationary probability that the system is in a state belonging to $S_i.$

If the system is in the set $S_0,$ then the probability that, at next step, the system will be in this set is $1-N\varepsilon+o(\varepsilon),$  the probability that the system will be in the set $S_1$ is  $N\varepsilon+o(\varepsilon),$ and the probability that the system will be in a state not belonging to
$S_0\cup S_1$ is $o(\varepsilon),$ $\varepsilon\to 0.$ The probability that, from the set $S_i,$ the system passes to the set $S_{i-1},$ $i=2,\dots,[N/2],$ is equal to $1-o(\sqrt(\varepsilon)).$ From this, it follows to  
$$P_0=1-N\varepsilon+o(\varepsilon),\eqno(14)$$  
$$P_1=N\varepsilon+o(\varepsilon),\eqno(15)$$
$$P_i=o(\varepsilon).\eqno(16)$$

Using (14)--(16), we get
$$v=P_0+\frac{(N-1)P_1}{N}+o(\varepsilon)=1-N\varepsilon+o(\varepsilon).$$ 
Theorem 2 has been proved.

\section{Conclusion}

For a deterministic binary chain of contours, we have obtained a competition resolutuon rule such that this rule is optimal in the following sense. From any initial state, the system results in a state of free movement over the minimum time interval. A stochastic version of the system is considered. In this version, the probability that an attempt of a particle to move is $1-\varepsilon.$ We have proved that the average velocity of particles tends to 1 as $\varepsilon\to 0.$ Under the assumption that the number of contours equals 3, we study the asymptotic behavior  
of the system with the left-priority resolution rule as $\varepsilon \to 0.$ For this rule, the average velocity of particles does not tend to 1 as $\varepsilon\to 0.$

\bibliographystyle{unsrt}  


\begin{thebibliography}{1}

\bibitem{SchrShN}   
Schreckenberg M., Shadshneider M., Nagel K., Ito N. Discrete stochastic models for traffic flow. Physical Review, 1995, vol.~51, no.~4, 2939. 

\bibitem{Blank}   
 Blank M.L. Exact analysis of dynamical systems arising
in models of traffic flow. Russian Mathematical
Surveys, 2000; 55:3, 562--563. \newblock {DOI: 10.4213/rm95}

\bibitem{BelFe}
 Belitsky V., Ferrari P.A. Invariant measures and
convergence properties for cellular automation 184 and related processes. 
J. Stat. Phys. 2005; vol.~118, no.~3/4, pp.~589~-- 523. DOI:10.1007/s10955-004-8822-4

\bibitem{GrGr}
 Gray L., Griffeath D. The ergodic theory of traffic jams.
J.~Stat.~Phys. 2001. V.~105. no.~3/4. P.~413~-- 452. DOI 10.1023/A:1012202706850

\bibitem{KanaiNiTo}
 Kanai M., Nishinary K., Tokihiro T. Exact solution and asymptotic behavior of the asymmetric behavior of the asymetric simple exclusion process on a ring. arXiv.0905.2795v1 [cond-mat-stat-mech] 18 May 2009.

\bibitem{Kanai}
 Kanai M. Two-line traffic-flow model with an exact steady-state solution, Physical Review, 2010, vol.~82, no.~6, 066107.   

\bibitem{Blank}
 Blank M. Metric properties of discrete time exclusion type processes in continuum. J. Stat. Phys. 2010, vol.~140,  no.~1.~-- pp.~170--197.

\bibitem{YT-2018}
 Yashina M.V., Tatashev A.G. Traffic model based on synchronous and asynchronous exclusion processes. Math. Method. Appl. Sci., 2018, vol.~43, issue~14, pp.~8136--8146.\\ DOI:~10.1002/mma.6237   

\bibitem{ Wolfram}
 Wolfram S. Statistical mechanics of cellular automata. Rev. Mod. Mod. Phys. 1983, vol. 55, pp. 601~-- 644.\\
DOI: Rev.Mod.Phys.55.601

\bibitem{Spitzer}
 Spitzer F. Interaction of Markov processes. {\it Advances in Mathematics,} vol. 5, 1970, pp. 246~-- 290.\\ DOI: 10.1016/0001-8708(70)90034-4

\bibitem{Biham}
Biham O., Middleton A.A., Levine D. Self-organization and a dynamic transition in traffic-flow models.  Phys. Rev.~A. 
The American Physical Society, 1992, vol.~46, no.~10. R6124--R6127.

\bibitem{Angel}
 Angel O., Horloyd A.E., Martin J.B. The jammed phase of the Biham-Middelton-Levine traffic flow model // Electronic Communication in Probability. 2005. V.~10. P.~167~-- 178.

\bibitem{BY}
Buslaev A.P., Yashina M.V. On holonomic mathematical $F$‐pendulum. Math. Method. Appl. Sci., 2016, vol.~39, issue~16, pp.~ 4820--4828. DOI: 10.1002/mma.3810

\bibitem{BBKY}
 Bugaev A.S., Buslaev A.P., Kozlov V.V., Yashina M.V. Distributed problems of monitoring and modern approaches to traffic modeling, p. 6, 14th International IEEE Conference on Intelligent Transportation Systems (ITSC 2011),
Washington, USA, 5-7.10.2011.  DOI: 10.1109/ITSC.20116082805. (2011) 477--481. 

\bibitem{KBT-2013}
 Kozlov V.V., Buslaev A.P., Tatashev A.G. On synergy of totally connected flows on chainmails /  Proceed. of International Conference of CMMSE, 24.06-27.06 2013, v.3,
pp.~861--874.

\bibitem{KBT-2015-CM}
 Kozlov V.V., Buslaev A.P., Tatashev A.G. Monotonic walks on a necklace and a coloured dynamic vector. International Journal of  Computer Mathematics, 2015, vol.~92, no.~9, pp. 910--920.  DOI: 1080/00207160.2014.915964

\bibitem{KBT-2015-JCAM}
 Kozlov V.V., Buslaev A.P., Tatashev A.G. A dynamic communication system on a network. Journal of  Computational and Applied Mathemstics, 2015, pp.~247--261.\\ 
DOI:~10.1016/j.cam.2014.07.026 

\bibitem{TY}
 Tatashev A.G., Yashina M.V. Spectrum of Elementary Cellular Automata and Closed Chains of Contours // Machines 2019, 7(2), 28.

\bibitem{BTFY}
  Buslaev A.P., Fomina M.Yu., Tatashev A.G., Yashina M.V. On discrete flow networks model spectra: statements, simulation, hypotheses. Journal of Physics: Conference Series, 1053 (2018) 012034, pp.~1--7. 
DOI:~10.1088/1742/6596/1053/1/012034

\bibitem{Blank}
 Bugaev A.S., Tatashev A.G.,  Yashina  M.V. Spectrum of a continuous closed symmetric chain with an arbitrary number of contours. Mathematical Models and Computer Simulations, 2021, 13(6), 1014-1027.

\bibitem{YT}
 Yashina M.V., Tatashev A.G. Spectral cycles and average velocity of clusters in discrete two-contours system with two nodes. Math. Method. Appl. Sci., 2020, vol.~43, issue7, pp.~4303--4316.
 DOI: 10.1102/mma6194

\bibitem{YTF}
  Yashina M.V., Tatashev A.G.,  Fomina M.Y. Optimization of velocity mode in Buslaev two-contour networks via competition resolution rules. International Journal of Interactive Mobile Technologies, 2020, vol.~14, no.~10, pp.~61--73.\\
DOI: 103991/ijim.v14i10.14641 

\bibitem{MTY}
 Myshkis P.A., Tatashev A.G.,  Yashina M.V. Cluster motion in a two-contour system with priority rule for conflict resolution.
Journal of Computer and Systems Sciences International, 2020, vol.~59, No.~3. pp.~311--321.

\bibitem{Kemeny}
 Kemeny J.G., Snell J.L  Finite Markov Chains. Springer-Verlag, New York, Berlin, 
 Heidelberg, Tokyo 1976.
 







\end{thebibliography}

\end{document}